\documentclass{amsart}
\usepackage{amsmath,amssymb,amsfonts,amsthm}
\usepackage{geometry}
\usepackage{hyperref}
\usepackage{xcolor}
\usepackage{enumitem}
\usepackage{graphicx}
\usepackage{subcaption}
\usepackage{changepage}
\usepackage{relsize}
\usepackage{comment}
\usepackage{caption}
\usepackage{subcaption}
\usepackage{array}
\usepackage{booktabs}
\usepackage{diagbox}
\usepackage{bm}
\usepackage{algorithm}
\usepackage{algpseudocode}
\usepackage{graphicx}

\theoremstyle{definition}
\newtheorem{theorem}{Theorem}[section]

\newtheorem{definition}[theorem]{Definition}

\theoremstyle{remark}

\numberwithin{equation}{section}

\makeatletter
\@for\@tempa:={A,B,C,D,E,F,G,H,I,J,K,L,M,N,O,P,Q,R,S,T,U,V,W,X,Y,Z}\do{%
  \expandafter\edef\csname b\@tempa\endcsname{\noexpand\mathbb{\@tempa}}%
  \expandafter\edef\csname c\@tempa\endcsname{\noexpand\mathcal{\@tempa}}%
}
\makeatother

\title[Kernel Regression Estimators on the Circle]
{Nonparametric Regression on the Circle with Applications to Climate and Environmental Science}

\author{Athanasios G. Georgiadis}
\address{School of Computer Science and Statistics,
Trinity College of Dublin}
\email{georgiaa@tcd.ie}
\author{Andrew P. Percival}
\address{School of Computer Science and Statistics,
Trinity College of Dublin}
\email{aperciva@tcd.ie}

\begin{document}

\begin{abstract}
Two important indicators of climate change are the decline in Arctic sea ice area and the increase in average temperatures globally. Both sea ice and temperature exhibit seasonal patterns, making circular regression a natural approach for their analysis. In this paper, we construct computationally tractable Nadaraya–Watson type kernel regression estimators on the unit circle, using rate-optimal kernel density estimators. We introduce a leave-one-out cross-validation procedure for the data-driven selection of the smoothness parameter and evaluate the performance of the proposed estimators through a simulation study covering 48 scenarios, including comparisons with existing local polynomial estimators. We then apply our methodology to historical Arctic sea ice area data and daily temperature observations in Stockholm. The fitted regression functions and their derivatives or finite differences allow us to investigate changes in both the levels and seasonal dynamics of these environmental indicators. In particular, we observe a decline in fitted sea ice area and substantial changes in its seasonal rates of growth and retreat, together with higher fitted temperatures across most of the year in Stockholm during the more recent period in comparison to historical data.
\end{abstract}

\maketitle

\section{Introduction} 
\label{sec: intro}
A key indicator of ongoing climate change is the decrease in Arctic sea ice (\cite{DriversSeaIce}). Since 1979 the annual mean Arctic sea ice extent has experienced a linear decrease with some climate model projections predicting that the Arctic Ocean will experience a summer without ice by 2035 (\cite{DK21}). The decline of sea ice levels has severe ecological consequences, see \cite{Post13} and \cite{SL16}.

Climate change is also associated with an increase in average global temperatures. Indeed, the work of \cite{Forster26} tells us that relative to measurements from 1859-1900, human-induced global warming had increased temperatures by $1.37^\circ$C by 2025. Therefore, at a given point of the Earth, we expect present-day temperatures to have a different distribution to past data. An increase in temperature has direct effects on human health, leading to an increase in heat-related death (\cite{Relvas26}), and further negative consequences on agriculture (\cite{Janni24}), biodiversity (\cite{Hab22}), and the spread of disease (\cite{Anikeeva24}).

We propose the use of nonparametric regression estimators to analyse changes in the annual patterns of sea ice and daily temperatures over time. Since both phenomena exhibit annual periodicity, it is natural to consider regression estimators whose domain is the unit circle $\mathbb{S}^1$, identifying the beginning and the end of each year. Indeed, in the Arctic, sea ice experiences a period of thawing, between March and September, and a period of refreezing for the remainder of the year. Thus, we can represent the time of year as a point on the unit circle, taking sea ice area as the real-valued response variable. Similarly, at a given point on Earth, daily temperatures tend to increase during summer and decrease during winter. Therefore, we can again represent the time of year on the unit circle, with daily temperature as the response variable.

Moreover, regression estimators with circular predictors have potential applications in several scientific fields. Examples of circular or periodic variables include lunar phases, wind and wave directions, seasonal weather patterns, diurnal temperature variations, circadian rhythms, hormonal cycles, and seasonal patterns of influenza infections. In such settings, one may wish to investigate the relationship between a circular variable $X$ and a real-valued response $Y$. More precisely, given $n$ i.i.d. observations $(X_1,Y_1),\dots,(X_n,Y_n)$ of a random pair $(X,Y)$ taking values in $\mathbb{S}^1\times\mathbb{R}$, with $\mathbb{E}|Y|<\infty$, our aim is to estimate the regression function
\begin{equation*}
R(x)=\mathbb{E}(Y\mid X=x),\qquad x\in\mathbb{S}^1.
\end{equation*}

A nonparametric solution to this problem in the case where $X$ and $Y$ are both real-valued was separately proposed by \cite{Nad64} and \cite{Wat64}. Such Nadaraya-Watson (NW) estimators of a regression function use `kernels' from the density estimation problem. That is, such kernels are used as weights for a locally weighted average of observations of the random pair $(X,Y)$.


An extension of the NW estimator to the case where the predictor takes values on a closed Riemannian manifold was developed by \cite{Pell06}. The construction uses kernels introduced in \cite{Pell05} for density estimation, which depend on the local geometry of the manifold. Pelletier establishes pointwise mean squared error bounds with the same convergence rate as the corresponding NW estimators on Euclidean spaces. More recently, \cite{Bouz24} developed kernel-based conditional U-statistic estimators for predictors and responses taking values in general metric spaces, using density estimation kernels from \cite{CGKPP20}, \allowbreak \cite{CGP22}, \allowbreak \cite{CGW25}. This framework includes regression on the unit circle with real-valued responses, as considered in the present study.

An important and broadly used approach to nonparametric regression on spheres was developed by \allowbreak \cite{DiMarzio}, who introduced local polynomial estimators for predictors taking values on $\mathbb{S}^d$, $d\geq 1$. Their construction uses von Mises–Fisher type kernels, previously studied by \cite{Hall87} for kernel density estimation on spheres. The authors also investigate cross-validation for smoothing parameter selection and evaluate their estimators through simulations and real-data applications in the case of polynomial degrees $p=0$ and $p=1$. Such estimators are of particular interest to our study and are used for comparison in our simulation study. For $p=0$, the local polynomial estimator reduces to a NW-type estimator based on von Mises–Fisher kernels.

An interesting and related line of research, involving Euclidean predictors and circular responses, has been developed in \cite{Vila22,Vila26}. These works introduce goodness-of-fit tests for parametric circular regression models, using nonparametric estimators for comparison. Their methodology is supported by simulation studies, while \cite{Vila26} also presents an environmental application involving wave directions in the Adriatic Sea. Further nonparametric approaches to regression with circular or spherical responses can be found in \cite{DM13,Vila21a,Vila21b}.

Further related developments in nonparametric regression and conditional estimation can be found in
\cite{AlonsoPena21,Aswani11,BouzebdaNezzal24}, \allowbreak

\cite{BouzebdaTaachouche24,Ferraty07,MeilanVila24}, and the references therein.
\vspace{0.2cm}

The \textit{aim of this paper} is to develop computationally tractable Nadaraya–Watson type regression estimators on the unit circle and apply them to the analysis of historical Arctic sea ice area and daily temperature data. Our construction builds on the rate-optimal kernel density estimators developed in our recent work \cite{GP26}, where infinite spectral expansions were replaced by finite sums while preserving the optimal convergence rate for density estimation. This provides a natural basis for constructing regression estimators that can be efficiently implemented in practice. An important component of our methodology is the data-driven selection of the smoothness parameter through leave-one-out cross-validation, allowing the estimators to be fitted without requiring a predetermined smoothness parameter.

We investigate the performance of the proposed estimators through an extensive simulation study covering $48$ scenarios, including comparisons with the well-established local polynomial estimators of \cite{DiMarzio}, demonstrating the effectiveness and competitiveness of our approach. We then apply our methodology to Arctic sea ice area data and daily temperatures in Stockholm, considering different historical periods. Our objective is to identify and quantify changes in the annual patterns of these environmental indicators. Beyond estimating the regression functions, our methodology allows us to compare fitted curves across historical periods and investigate their \textit{seasonal dynamics through derivatives and finite differences}. This provides additional information about the magnitude and timing of seasonal changes, beyond differences in the levels of environmental indicators. Our results indicate substantial changes in the seasonal rates of Arctic sea ice retreat and refreezing. For Stockholm, the fitted temperatures during 2010–2025 are higher than those during 1940–1955 for approximately $97\%$ of the year. These findings illustrate the potential of nonparametric regression on the circle as a useful statistical tool for analysing long-term changes in environmental data with an annual periodic structure.
\vspace{0.2cm}

The remainder of this paper is organised as follows. Section \ref{sec: data} presents the Arctic sea ice area and daily temperature datasets, together with the necessary pre-processing steps. Section \ref{sec: method} introduces the statistical framework, reviews the relevant developments in kernel density estimation, and constructs the proposed kernel regression estimators, with their implementation summarised in Algorithm \ref{alg: KRE}. Section \ref{sec: simulations} evaluates the performance of the estimators across $48$ simulation scenarios and compares them with the local polynomial estimators of \cite{DiMarzio}. Section \ref{sec: results} presents and discusses the results of the environmental applications. Finally, Section \ref{sec: conclusion} provides the conclusions, while Section \ref{sec: data avail} contains the data and code availability statement.

\section{Data}
\label{sec: data}
\subsection{Dataset: Sea Ice Area in the Northern Hemisphere}
The interesting recent work of \cite{SL16} identifies a relationship between the sea ice concentration in Baffin Bay, Greenland, and the population of polar bears living in the region. Within this work, they plot the mean daily sea ice area for each decade alongside the the day of the year. They use an average curve in order to do this, giving rough estimates that are not differentiable. We are motivated by their work to use the same data and fit our kernel regression estimator (KRE) to it as an illustration of a potential use of such machinery. 

\begin{figure}
    \centering
    \includegraphics[width=1\linewidth]{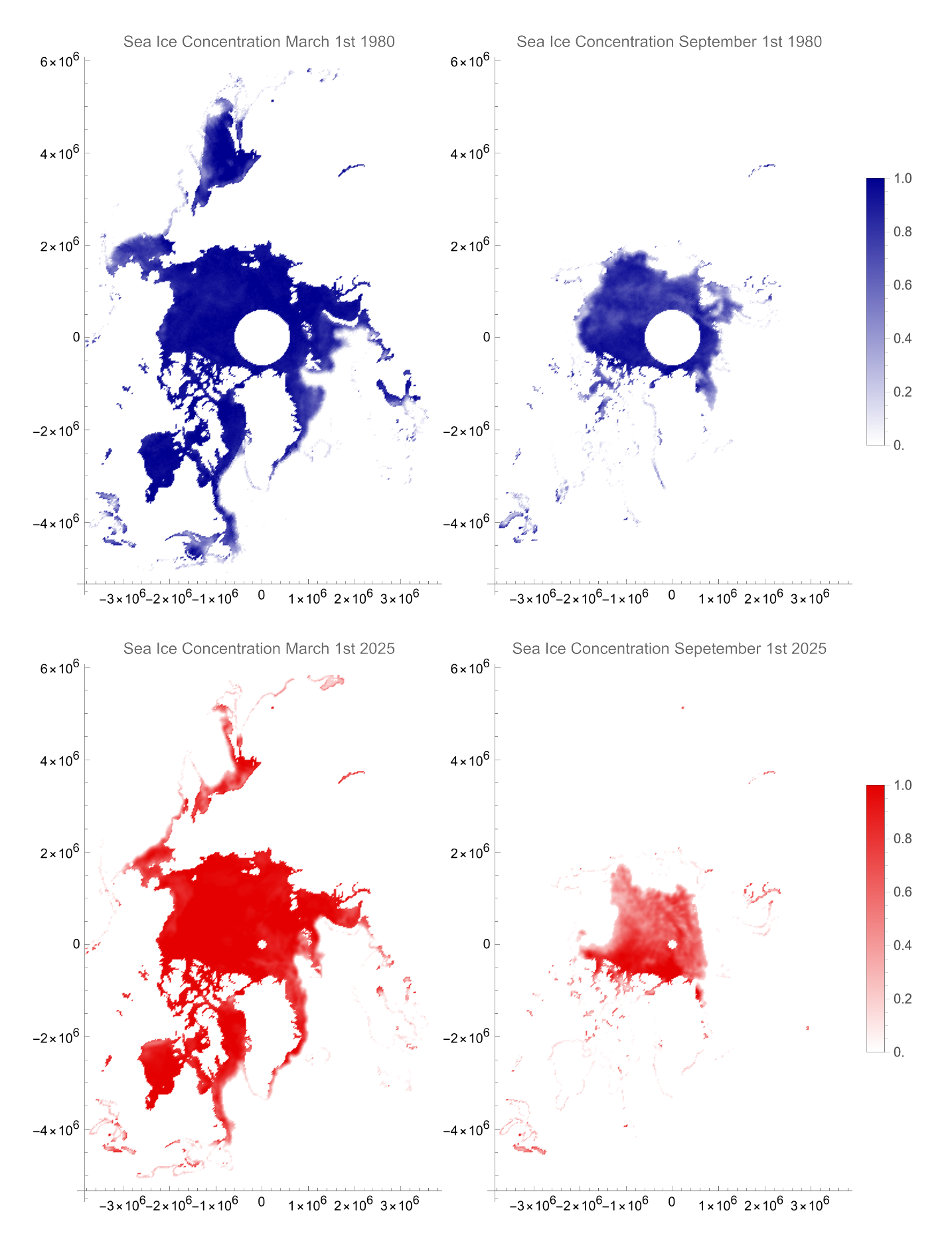}
    \caption{Plots of the sea ice concentration centred on the North Pole for the first of March and September in both 1980 and 2025. The axes give the distance (in metres) from the North Pole using the NSIDC polar stereographic projection. }
    \label{fig:northern hemisphere plot}
\end{figure}

Sea ice area in the northern hemisphere is constantly changing, but over the course of a year can be broken into a season of thawing, from March to September, and a season of freezing for the remainder of the year. In this way, the day of the year can be seen as the independent variable, taking periodic (i.e. circular) values over the course of the year, and the total sea ice area is our response variable. In Figure \ref{fig:northern hemisphere plot} we plot the sea ice concentrations centred on the North Pole for the dates of the first of March and first of September for both 1980 and 2025. One may see that the amount of sea ice in the 2025 plots is less than for the corresponding dates in 1980. This motivates us to construct a KRE for each decade of sea ice data in order to properly interpret this change. 

The data used is the Sea Ice Concentrations gathered from Nimbus-7 SMMR and DMSP SSM/I-SSMIS Passive Microwave Data, freely available through the National Snow and Ice Data Center (NSIDC), \url{https://nsidc.org/data/nsidc-0051/versions/2}. The data is provided in NSIDC's polar stereographic projection at a grid cell size of 25 $\times$ 25 km. We elect to use the entire dataset, not just the subset within Baffin Bay, which covers from $30.98^\circ$N northwards to between $84.5^\circ$N and $89^\circ$N, depending how recent the observation is. This `polar hole' is based on the limitations of satellite observations and has become smaller as technologies improve. This hole can also be seen in the centre of each plot in Figure \ref{fig:northern hemisphere plot}.

\subsubsection{Data pre-processing}

The dataset's response variable is the ice concentration, not area, given as a number between $0$ and $1$. The dataset also contains values greater than $1$. A value of $1.004$ corresponds to missing values from within the polar hole. Values greater than $1.004$ correspond to other various types of missing or irrelevant data (such as coasts and landmasses). We elect to set all values greater than $1$ to zero. In particular, for the polar hole, the oldest dataset has is missing observations from $84.5^\circ$ northwards. Therefore, we set all observations of sea ice concentration above $84.5^\circ$N to zero. This is to ensure that the sea ice area is calculated over the same geographical region across all historical periods. So, to be more precise, our study concerns sea ice area between $30.98^\circ$N and $84.5^\circ$N. 

We then set all concentrations less than $0.15$ to zero. This is because a satellite reading of less than $0.15$ becomes unreliable given noise e.g. from cloud-cover. This follows the convention of \cite{Zwally83}, which is noted in \cite{Meier26} to give an average agreement between observed sea ice area and satellite records.

These concentrations are then multiplied by the area of each grid cell. Note that these are not all $25\times25=625$ km$^2$ due to the curvature of the Earth. We read the true areas from the file \texttt{NSIDC0771\_CellArea\_PS\_N25km\_v1.1.nc}, also freely available. The sum of each pixel's concentration multiplied by its area gives the observation of sea ice area for each day. 

The dataset spans from the 26th of October 1978 to the 31st of December 2025. From the 26th of October 1978 until the 19th of August 1987 (inclusive), there are only observations every second day. From the 20th of August 1987 until the 2nd of December 1987 (inclusive) there is an observation every day. There are no available observations from the 3rd of December 1987 until the 12th of January 1988 (inclusive). From the 13th of January 1988 until the 31st of December 2025, there is an observation on each day. The data is converted into angles in radians using Algorithm \ref{alg: day to angles}. In this way, we have observations $\theta_1,\dots,\theta_n\in(0,2\pi]$ which are suitable as inputs into the KREs given later in Algorithm \ref{alg: KRE}.

\begin{algorithm}
\caption{Converting dates to observations on $\mathbb{S}^1$}\label{alg: day to angles}
\begin{algorithmic}[1]
    \Require Vector $T$ of dates of observations
\State $n\gets \text{length}(T)$
\State Create vector $\theta$ of length $n$
\State Define function $d(t):=\text{`day of year of } t \text{'}$ for any date $t$
\State Define function $\ell(t):=\begin{cases}
    1,\quad \text{if }t \text{ is in a leap year,}\\
    0, \quad \text{otherwise},
\end{cases}$ for any date $t$
\State $i\gets1$
\While{$i \leq n$}
\State $\theta_i\gets \dfrac{2\pi d(T_i)}{365+\ell(T_i)}$
\State $i\gets i+1$
\EndWhile

\State \Return Vector $\theta$
   
\end{algorithmic}
\end{algorithm}

Following the work of \cite{SL16}, we do not use the entire dataset in our KRE and instead elect to split the data into five groups; from 1978 to 1987 ($n=1715$), from 1988 to 1997 ($n=3639$), from 1998 to 2007 ($n=3651$), from 2008 to 2017 ($n=3653$), and from 2018 to 2025 ($n=2922$). This allows us to create a KRE for each period in order to compare to change in the thawing-freezing cycle across decades.

\subsection{Dataset:  Daily Mean Temperature in Stockholm}
Our second case study concerns estimating the daily mean temperature across the calendar year at a particular point on the Earth, in this case Stockholm, Sweden. However, due to the effects of climate change one cannot assume that the temperature at a particular point on the Earth at a particular time in the Earth's orbit has the same expected value now as it did eighty years ago. Indeed, in this application, we will create two KREs based on observations of daily mean temperature in Stockholm from two different periods: 1940-1955, and 2010-2025. In order to demonstrate climate change's effect of increasing average global temperatures, we will show that each dataset produces quite different KREs. 

This data is freely available from the Open-Meteo Historical Weather API, \url{https://open-meteo.com/en/docs/historical-weather-api}. The dataset in question consists of observations of the daily mean temperature (recorded in degrees Celsius) in Stockholm from the 1st of January 1940 to the present. The date data is converted to angles in radians via Algorithm \ref{alg: day to angles}. This leaves us with $n=5844$ observations in each of the 1940-1955 and 2010-2025 periods.

\section{Methodology}
\label{sec: method}

\subsection{Theoretical background}
\label{sec: geom background}
In this study, the pair of random variables $(X,Y)$ take values on $\mathbb{S}^1\times\mathbb{R}$ where $\mathbb{S}^1=\{x\in \mathbb{R}^2 \,:\, \|x\|=1\}$ is the unit circle. A given point $x=(x_1,x_2)\in\mathbb{S}^1$ has a polar angle representation $\theta\in(-\pi,\pi]$ given by
\begin{equation}
    \label{eq: Cart to theta}
    \theta=\begin{cases}
        \arccos x_1,\quad &\text{if }x_2\geq 0,\\
        -\arccos x_1,\quad &\text{if }x_2<0.
    \end{cases}
\end{equation}
Conversely, for a given angle $\theta\in(-\pi,\pi]$, we have a unique point on the circle given by
\begin{equation}
    (x_1,x_2)=(\cos\theta,\sin\theta).
\end{equation}

Therefore, let us fix the following notational convention throughout the remainder of this work that for a set of observations $X_1,\dots,X_n\in\mathbb{S}^1$, for some $n\in\mathbb{N}$, let $\theta_i\in(-\pi,\pi]$ be the polar angle of $X_i$ given by \eqref{eq: Cart to theta}, for $1\leq i\leq n$.

By identifying the endpoints of $(-\pi,\pi]$, random variables on the circle can be represented by angular variables, and the methodology developed in this study extends naturally to random variables with any period. Indeed, consider a design $T_1,\dots,T_n$ taking values in $(a,b]$ with density $f_T$ on this interval, whose periodic extension to $\mathbb{R}$ has period $b-a$. Then, by Algorithm \ref{alg: T to angles}, we can transform these observations into angular variables $\theta_1,\dots,\theta_n$ taking values in $(-\pi,\pi]$, whose density admits a $2\pi$-periodic extension to $\mathbb{R}$. These can therefore be regarded as observations on $\mathbb{S}^1$.

\begin{algorithm}
\caption{Converting periodic data to observations on $\mathbb{S}^1$}\label{alg: T to angles}
\begin{algorithmic}[1]
    \Require Vector $T$ of observations periodic on some $(a,b]$
\State $n\gets \text{length}(T)$
\State Create vector $\theta$ of length $n$
\State $i\gets1$
\While{$i \leq n$}
\State $\theta_i\gets \dfrac{2\pi}{b-a}\cdot (T_i-a)-\pi$
\State $i\gets i+1$
\EndWhile

\State \Return Vector $\theta$
   
\end{algorithmic}
\end{algorithm}

Therefore, the study of random variables on the unit circle is well motivated, as it includes not only inherently circular variables, such as wind direction, but also periodic variables defined over a given interval, such as the time of day, the calendar year, or the orbital position of a celestial body.
\vspace{0.2cm}

In nonparametric estimation, the usual approach is to assume that the unknown function belongs to some smoothness class. One such class is the H{\"o}lder space of order $s>0$, which is also considered in \cite{GP26}. Note that throughout this study, $\lfloor s\rfloor$ represents the greatest integer \textit{strictly} less than $s$. Similarly, $\lceil s\rceil$ is the smallest integer strictly greater than $s$.
\begin{definition}
    Let $s>0$ and $f:\mathbb{S}^1\to\mathbb{R}$ be a function. We say that $f$ belongs to the H{\"o}lder space $\dot{\mathcal{H}}^s$, if \begin{equation}
        \|f\|_{\dot{\mathcal{H}}^s}:=\sup_{x\neq y}\dfrac{|(-\Delta)^{\lfloor s \rfloor/2}f(x)-(-\Delta)^{\lfloor s\rfloor/2}f(y)|}{\rho(x,y)^{s-\lfloor s\rfloor}}<\infty,
    \end{equation}
    where $\rho(x,y):=\arccos(\langle x,y\rangle)$ is the distance on $\mathbb{S}^1$, and the Laplacian in polar coordinates is $\Delta f:=d^2f/d\theta^2$.
\end{definition}
Further details on the fractional powers of the Laplace-Beltrami operator can be found in \cite{ReedSimon}.

\subsection{Kernel density estimation on the circle}
\label{sec: KDEs}

Our regression estimator is constructed using kernels developed for density estimation on the unit circle. We briefly recall the relevant results from \cite{CGKPP20,CGW25,GP26}.

Consider i.i.d. observations $X_1,\dots,X_n$ on $\mathbb{S}^1$, with unknown density $f_X$. For an appropriate even function $k:\mathbb{R}\to\mathbb{R}$ and bandwidth $h>0$, the corresponding KDE is given by
\begin{equation}
    \widehat{f}_{n,h}(\theta)
    =\frac{1}{2\pi n}\sum_{i=1}^n
    \Big(1+2\sum_{\nu=1}^{\infty}
    k(h\nu)\cos(\nu(\theta-\theta_i))\Big),
    \quad \theta\in(-\pi,\pi].
\end{equation}

The infinite spectral expansion is computationally inefficient.  In \cite{GP26} we show that it can be replaced by a finite sum while preserving the optimal convergence rate. In particular, for $s>0$, let $h=n^{-1/(2s+1)}$ and
\begin{equation}
\label{eq: truncation}
    N_s:=
    \Big\lfloor
    \Big(\frac{n^{(s+r)/(2s+1)}}{\pi(r-1)}\Big)^{1/(r-1)}
    \Big\rfloor+1,
\end{equation}
where $r$ is determined by the decay properties of $k$. The resulting estimator is
\begin{equation}
\label{eq: KDE truncated}
    \widehat{f}_s(\theta)
    :=\frac{1}{2\pi n}\sum_{i=1}^n
    \Big(1+2\sum_{\nu=1}^{N_s}
    k(h\nu)\cos(\nu(\theta-\theta_i))\Big).
\end{equation}

In particular, we consider the function
\begin{equation}
\label{eq: symbol}
    k(\lambda)=\frac{1}{1+|\lambda|^r},
    \qquad r=\lceil s\rceil+3.
\end{equation}
For this choice of $k$, the results of \cite{GP26} establish that, if $f_X\in\mathbb{L}^{\infty}(\mathbb{S}^1)\cap\dot{\mathcal{H}}^s(\mathbb{S}^1)$, then there exists some $C>0$ such that
\begin{equation}
    \sup_{\theta\in(-\pi,\pi]}
    \mathbb{E}\Big[
    \big(\widehat{f}_s(\theta)-f_X(\theta)\big)^2
    \Big]
    \leq C n^{-2s/(2s+1)}.
\end{equation}
Thus, the finite-order KDE is computationally tractable and retains the optimal convergence rate. For further details and the general assumptions on $k$, we refer to \cite{GP26}.

\subsection{Kernel regression estimators on the circle}
\label{sec: KREs}
With suitable kernels in hand, we will present the regression estimators on the circle. Firstly, let us illustrate the idea underlying this estimator.

Assume $X$ is a random variable with sample space $\mathcal{M}$, a metric measure space. Assume that the unknown density $f_X$ is strictly positive on $\mathcal{M}$. Furthermore, let $Y$ be a real-valued random variable associated with $X$. Suppose there exists a measurable \textit{regression function} $R:\mathcal{M}\to\mathbb{R}$ such that
\begin{equation}
\label{eq: regression function definition}
    R(x)=\mathbb{E}(Y|X=x),\quad \forall x\in \mathcal{M},
\end{equation}
Assume that there exists a real-valued random variable $\xi$ such that $\mathbb{E}(\xi|X)=0$ and $\text{Var}(\xi)<\infty$,
i.e. the regression model is given by 
\begin{equation}\label{eq:regression_model}
Y = R(X) + \xi.
\end{equation}

We observe $n$ independent realizations
$
(X_1,Y_1),\dots,(X_n,Y_n),
$
of the pair $(X,Y)$, each satisfying model \eqref{eq:regression_model}.
The nonparametric regression problem consists of finding an estimator $\widehat{R}_n$ for the unknown regression function $R$. The classic Nadaraya-Watson approach is the following. 

By (\ref{eq: regression function definition}) we re-express $R({x})$ as
\begin{equation}
\label{eq: regression with densities}
R({x})=\int_{\bR}yf_{Y|{X}}(y|{x})dy=\dfrac{1}{f_X(x)}\int_\mathbb{R} y f_{XY}(x,y)dy,
\end{equation}
where $f_{Y|{X}}(y|{x})$ is the conditional density of $Y$, given that $X={x}$, and $f_{XY}$ is their joint density. Now we just replace the above densities with appropriate kernel density estimators.

Assume that $\cK_h:\cM\times\cM\to\bR$, is a suitable kernel function and $h>0$ is the corresponding bandwidth. Denote by 
\begin{equation}\label{eq: KDE M}
\widehat{f}_n({x}):=\frac{1}{n}\sum_{i=1}^{n}\cK_h({X}_i,{x}),\qquad {x}\in\cM,
\end{equation}
the associated kernel density estimator for the unknown $f_{{X}}$. 

Let also $K:\bR\to\bR$ be an \textit{even} function with $\int_{\bR}K(u)du=1$ and $h'>0$ the corresponding bandwidth. An appropriate estimator for the joint density  $f_{XY}$ is as follows:
\begin{equation}
\widetilde{f}_n(x,y)=\frac{1}{n}\sum_{i=1}^{n}\cK_h({X}_i,{x})\frac{1}{h'}K\Big(\frac{y-Y_i}{h'}\Big),\qquad x\in\cM,\;y\in\bR.
\end{equation}
We substitute the two estimators and perform the integration in (\ref{eq: regression with densities}), by changing the variable and using that $\int K(u)=1$, and $\int uK(u)=0$, (because $uK(u)$ is an odd function) to receive the final regression estimator of Nadaraya-Watson type:
\begin{equation}
\label{eq: NW general}
\widehat{R}_n({x})=\frac{\sum_{i=1}^{n}Y_i\cK_h({X}_i,{x})}{\sum_{i=1}^{n}\cK_h({X}_i,{x})},\qquad x\in\cM.
\end{equation}
The construction above is not restricted to the unit circle. It can be applied to other metric measure spaces $\mathcal{M}$ for which suitable kernel density estimators of the form \eqref{eq: KDE M} are available, for example the interval, the sphere, the ball, and Riemannian manifolds (see \cite{CGP22}). Therefore, with the computationally efficient density estimator \eqref{eq: KDE truncated} on $\mathbb{S}^1$ at hand, we can present the following regression estimator for each $s>0$:

\begin{equation}
\label{eq: KRE on S1}
    \widehat{R}_{s}(\theta):=\dfrac{\sum_{i=1}^nY_i(1+2\sum_{\nu=1}^{N_s} k(h\nu)\cos(\nu(\theta-\theta_i))}{\sum_{i=1}^n(1+2\sum_{\nu=1}^{N_s} k(h\nu)\cos(\nu(\theta-\theta_i))},\quad \theta\in(-\pi,\pi],
\end{equation}
for i.i.d. observations $(X_1,Y_1),\dots,(X_n,Y_n)$, for some $n\in\mathbb{N}$, of the random pair $(X,Y)$ with sample space $\mathbb{S}^1\times\mathbb{R}$, where $h=n^{-1/(2s+1)}$, $k$ is given by (\ref{eq: symbol}). The proposed estimators provide a finite and therefore computationally efficient
realisation of the kernel-based regression framework considered in \cite{Bouz24}.

Since the bandwidth is determined by $h=n^{-1/(2s+1)}$, and the symbol $k$ and truncation level $N_s$ also depend on $s$, the KRE is fully determined by the smoothness parameter $s>0$. Thus, the data-driven selection of the smoothing parameter is transferred from the bandwidth $h$ to the smoothness parameter $s$. We select $s$ using leave-one-out cross-validation:

For a given $s>0$ and $i\in\{1,\dots,n\}$, let $\widehat{R}_{s,-i}$ denote the KRE as in \eqref{eq: KRE on S1} using a sample size of $n-1$, with the observation $(X_i,Y_i)$ excluded. Define the leave-one-out\footnote{We investigated alternative splits of training/test data (90:10, 80:20, 70:30, etc. and empirically found that none gave performance as good as leave-one-out)} cross-validator
\begin{equation}
\label{eq: CV for KRE}
    \text{CV}(s):=\dfrac{1}{n}\sum_{i=1}^n\big(Y_i-\widehat{R}_{s,-i}(\theta_i)\big)^2.
\end{equation}
Cross-validation is a well-established method for selecting smoothing parameters in nonparametric regression, motivated by its connection to the minimisation of prediction error and integrated squared error; see \cite{HM85}. Therefore, for a given dataset, we select the value of $s>0$ that minimises $\mathrm{CV}(s)$. The entire KRE construction is summarised in Algorithm \ref{alg: KRE}. Note that the Algorithm does not require the observations $\theta_i$ to lie in the interval $(-\pi,\pi]$. This is because the KRE, as a function of any given $\theta_i$ is $2\pi$-periodic. 

\begin{algorithm}
\caption{Kernel regression estimator construction}\label{alg: KRE}
\begin{algorithmic}[1]
    \Require Pair of vectors $(\theta,Y)$ of equal length 
    \Require Vector $S$ of candidates for the minimiser of the cross-validation function
\State $n\gets \text{length}(\theta)$
\State $M\gets \text{length}(S)$
\State Create vector $\text{CV}$ of length $M$
\State $ m\gets 1$
\While{$m\leq M$}
\State $s\gets S_m$
\State $h\gets (n-1)^{-1/(2s+1)}$
\State $r\gets\lceil s\rceil +3$
\State Define function $k(\lambda):=1/(1+|\lambda|^r)$
\State $N_s\gets \left\lfloor\left(\dfrac{ (n-1)^{(s+r)/(2s+1)}}{\pi(r-1)}\right)^{1/(r-1)}\right\rfloor+1$
\State $\text{CV}_m\gets \dfrac{1}{n}\mathlarger{\sum}_{i=1}^n\Bigg(Y_i-\dfrac{\sum_{j\neq i}Y_j(1+2\sum_{\nu=1}^{N_s} k(h\nu)\cos(\nu(\theta_i-\theta_j))}{\sum_{j\neq i}(1+2\sum_{\nu=1}^{N_s} k(h\nu)\cos(\nu(\theta_i-\theta_j))}\Bigg)^2$
\State $m\gets m+1$
\EndWhile
\State $\tilde{m}\gets \underset{1\leq m\leq M}{\text{argmin}} (\text{CV}_m)$
\State $s\gets S_{\tilde{m}}$
\State $h\gets n^{-1/(2s+1)}$
\State $r\gets\lceil s\rceil +3$
\State Define function $k(\lambda):=1/(1+|\lambda|^r)$
\State $N_s\gets \left\lfloor\left(\dfrac{ n^{(s+r)/(2s+1)}}{\pi(r-1)}\right)^{1/(r-1)}\right\rfloor+1$
\State Define function $\widehat{R}_s(t):=\dfrac{\sum_{i=1}^nY_i(1+2\sum_{\nu=1}^{N_s} k(h\nu)\cos(\nu(t-\theta_i))}{\sum_{i=1}^n(1+2\sum_{\nu=1}^{N_s} k(h\nu)\cos(\nu(t-\theta_i))}$
\State \Return Function $\widehat{R}_s(t)$
   
\end{algorithmic}
\end{algorithm}

\section{Simulations}
\label{sec: simulations}
Let us now proceed to test the performance of our KRE over several simulations. Firstly we shall examine the behaviour of the KRE across several scenarios wherein we vary the sample size and amount of noise in the response variable. We then proceed to compare the performance of this KRE against the local polynomial estimators of degree $p=0$ and $p=1$, from \cite{DiMarzio} across the same set of scenarios.

In this section we consider four $2\pi$-periodic regression functions $R^1,R^2,R^3,R^4:\bR\to\mathbb{R}$. For a given regression function $R^j$, $j\in\{1,2,3,4\}$, we simulate $n\in\{100,500,1000\}$ points from a uniform distribution on $\mathbb{S}^1$. This forms our observations $X_1,\dots,X_n$. In order to attain dependent data, for each $i\in\{1,\dots,n\}$ we sample $\epsilon_i\sim N(0,\sigma^2)$ for a fixed $\sigma\in\{0.5,1,2,3\}$. We then set $Y_i=R^j(\theta_i)+\epsilon_i$. In this way, we test regression estimators against 48 different scenarios; twelve for each regression function $R^j$. 

The regression functions are given as follows:
\begin{align*}
    &R^1(\theta)=5\cos\theta,\\ 
    &R^2(\theta)=5\sin2\theta,\\ 
    &R^3(\theta)=
1.18102 (3\sin^3(\theta)-\cos^2(\theta)+2\cos(3\theta))+0.14927,\\
&R^4(\theta)=-10.187+3.267\sum_{k=1}^4 a_k \exp[b_k \cos (c_k \theta)+d_k \sin (\theta)],
\end{align*}
where $\boldsymbol{a}=(0.3,1.5,-0.5,1)$, $\boldsymbol{b}=(1,-0.2,0,1)$, $\boldsymbol{c}=(5,1,1,1)$, and $\boldsymbol{d}=(0,-1,0.3,1)$. 
Each of these functions is $2\pi$-periodic on $\bR$ and therefore can be viewed as a function on $\mathbb{S}^1$.

Furthermore, note that each regression function achieves a global maximum and minimum at approximately $5$ and $-5$, respectively. This is to ensure that the added Gaussian noise to our simulated data has a comparable effect on the performance of each estimator.

In order to evaluate the performance of a given estimator, we use the mean integrated squared error (MISE). For a given function $R$ and its estimator $\widehat{R}$, define
\begin{equation}
    \text{MISE}(\widehat{R})\equiv\text{MISE}(\widehat{R},R):=\mathbb{E}\Big[\|\widehat{R}-R \|^2_2\Big]=\mathbb{E}\Big[\int_{-\pi}^\pi (\widehat{R}(\theta)-R(\theta))^2d\theta\Big].
\end{equation}
We estimate the MISE by taking the sample mean of the square of the $\mathbb{L}^2$-distance of 30 simulations for a fixed $n$ and $\sigma$. Importantly, when comparing the MISE of two or more estimators for a given $n$ and $\sigma$, we use the same thirty simulations. 

\subsection{Circular regression simulations}
\label{sec: regression sims subsec}
We now evaluate the performance of the KRE \eqref{eq: KRE on S1} by estimating its MISE as an estimator of the regression functions $R^j$, $j=1,2,3,4$, for the above regression model with sample size $n\in\{100,500,1000\}$ and noise level $\sigma\in\{0.5,1,2,3\}$.

In practice when selecting the value of $s>0$ that minimises $\text{CV}(s)$ as in \eqref{eq: CV for KRE}, we evaluate $\text{CV}(s)$ at $s\in\{0.4,0.45,0.5,\dots,3.5\}$. This upper and lower bound were determined heuristically. The upper is to ensure a finite search and the lower is because the cross-validation function become ill-behaved for very small values of $s$. We deemed this sufficient for our simulations and applications but for a given real-world dataset, further values of $s$ could be considered. 

\begin{figure}
    \centering
    \includegraphics[width=1\linewidth]{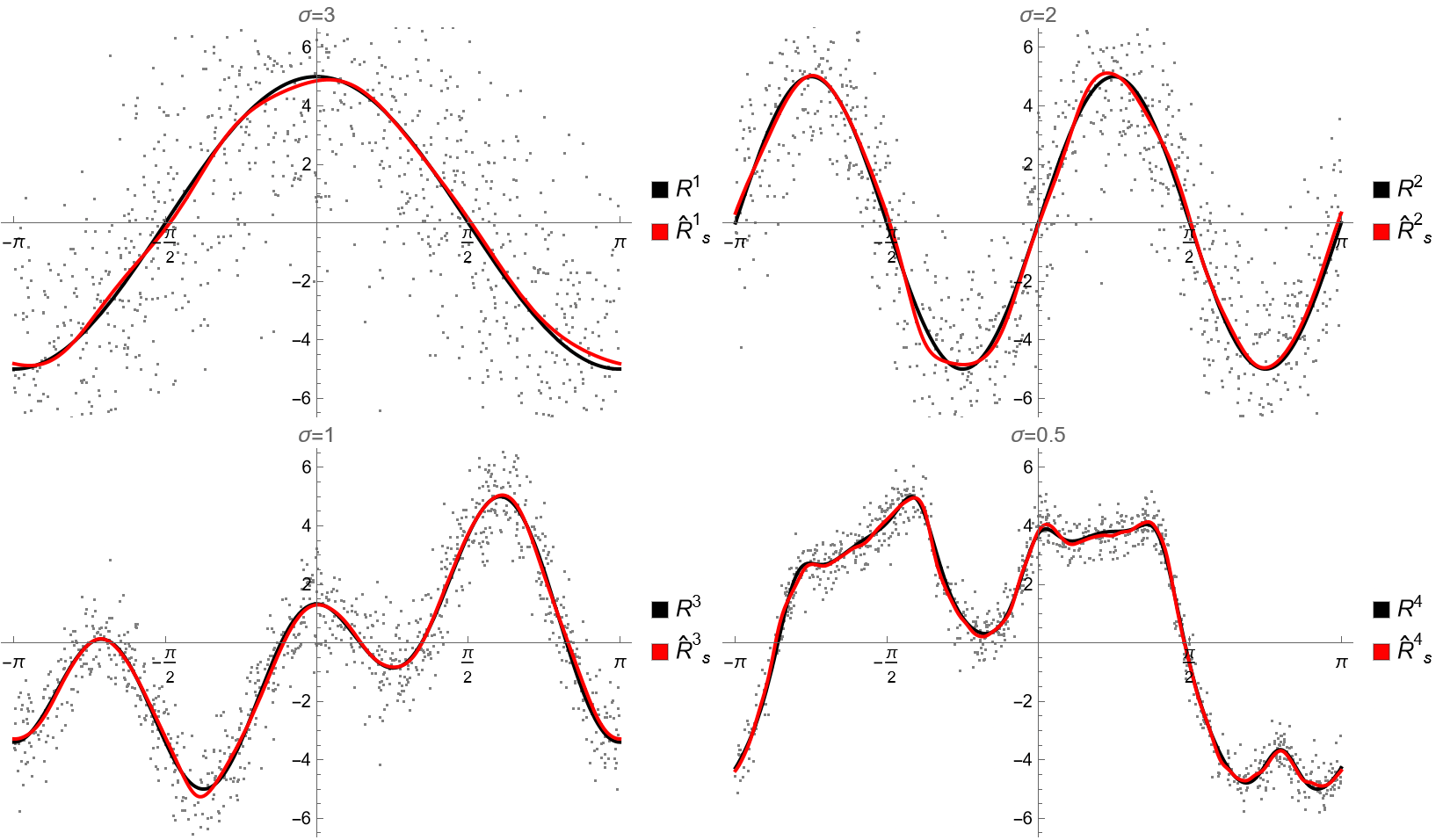}
    \caption{The true regression function with $n=1000$ points uniformly sampled with Gaussian noise (zero mean and $\sigma\in\{0.5,1,2,3\}$) and the corresponding finite order KREs $\widehat{R^j}_s$, $j\in\{1,2,3,4\}$, plotted in polar coordinates.}
    \label{fig: s1 KRE plots}
\end{figure}

In Figure \ref{fig: s1 KRE plots} we plot the regression functions $R^j$, $j\in\{1,2,3,4\}$, as well as the corresponding KRE $\widehat{R^j_s}$ from \eqref{eq: KRE on S1}, for a sample of $n=1000$. The standard deviation is taken to be $\sigma=3,\,2,\,1,\,0.5$ for $R^1,\,R^2,\,R^3,\,R^4$, respectively. 
We can see that in all cases, the KREs seem to provide a good fit, across noise levels. 

\begin{figure*}[]
    \centering
    \begin{subfigure}[t]{1\textwidth}
        \centering
        \includegraphics[width=1\linewidth]{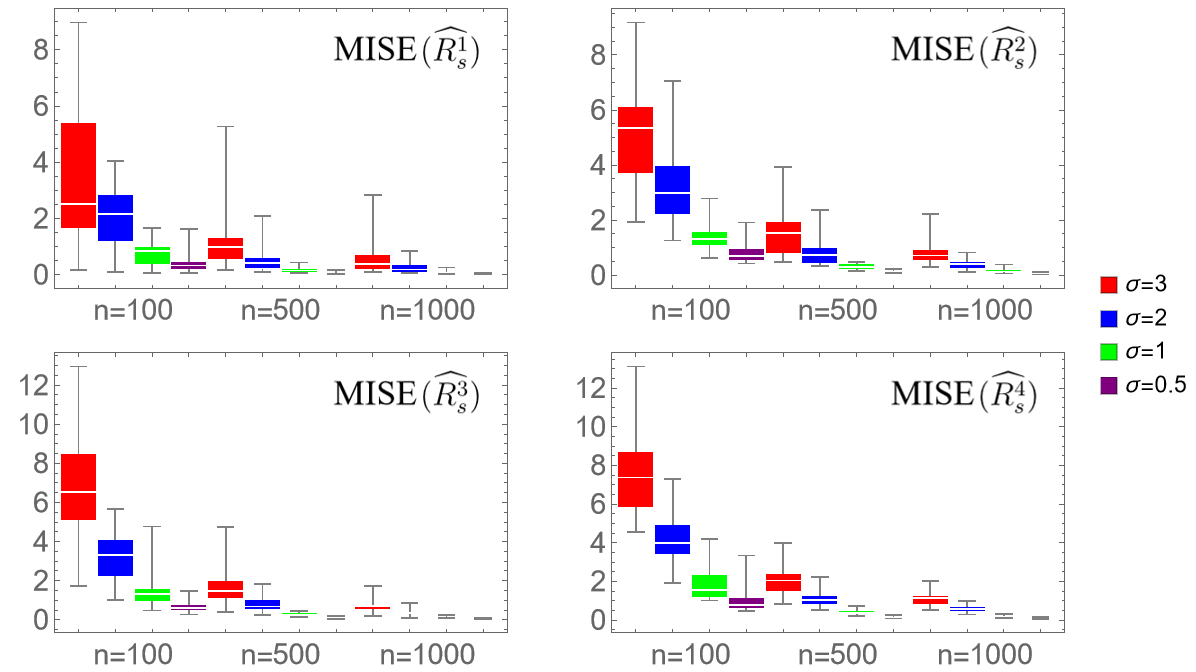}
        \caption{The MISE boxplots of the KRE for $R^j$, with $j\in\{1,2,3,4\}$, $n\in\{100,500,1000\}$, and $\sigma\in\{0.5,1,2,3\}$.}
    \end{subfigure}
    
        \begin{subfigure}[t]{1\textwidth}
        \centering
        \includegraphics[width=1\linewidth]{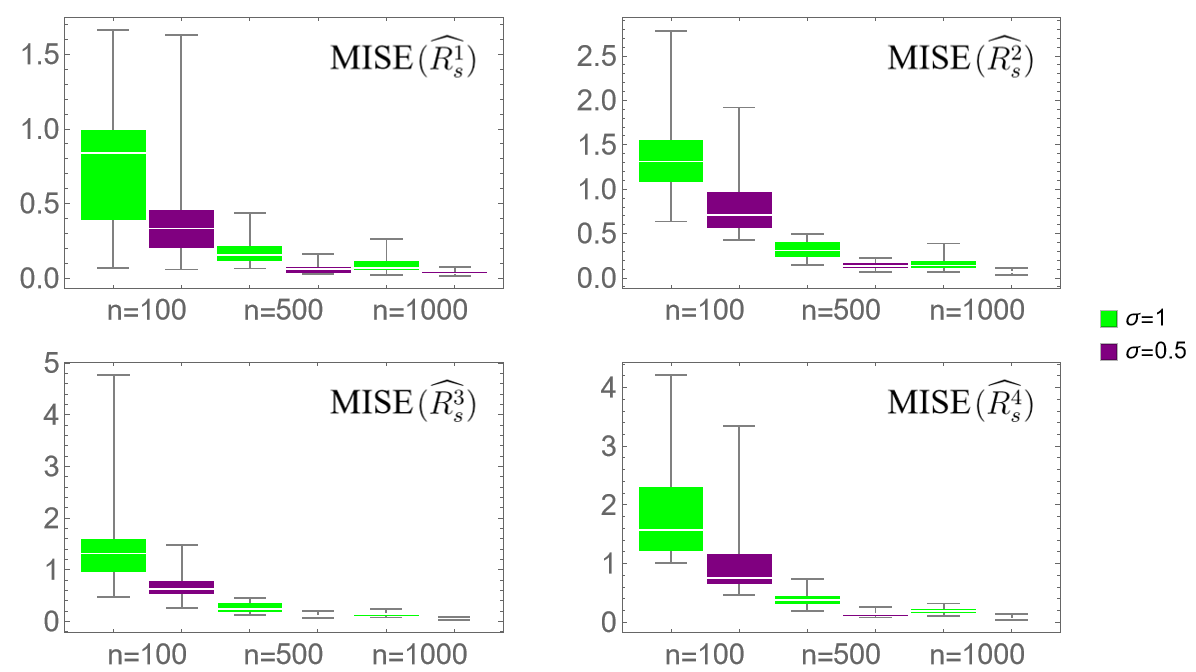}
        
        \caption{The MISE boxplots of the KRE for $R^j$, with $j\in\{1,2,3,4\}$, $n\in\{100,500,1000\}$, and $\sigma\in\{0.5,1\}$.}
    \end{subfigure}
    
    \caption{}
    \label{fig: s1 MISE plots}
\end{figure*}

We formally compare the performance of the KREs across these scenarios in Figure \ref{fig: s1 MISE plots}, where we give boxplots for the MISE of each estimator over $30$ simulations, for each regression function, sample size, and noise level. A subfigure of Figure \ref{fig: s1 MISE plots} gives a separate plot for the case of $\sigma\in\{0.5,1\}$ since the scale is much smaller.

In all cases, the result is similar: a higher sample size $n$, and lower standard deviation of noise $\sigma$, results in a better performance of the KRE, as expected. Note that for high sample sizes, the KREs are robust against large amounts of noise. 

\subsection{Comparison with local polynomial estimators}
\label{sec: LPEs}
We now compare our KRE with the local polynomial estimators (LPE) of \cite{DiMarzio}. These LPEs offer an alternative to the method presented in this study for regression estimation given a random pair $(X,Y)\in\mathbb{S}^1\times\mathbb{R}$ with i.i.d. sample $(X_1,Y_1),\dots,(X_n,Y_n)$, $n\in\mathbb{N}$. Di Marzio et al. present such LPEs in the case of the $d$-sphere, $d\geq 1$, for local polynomials of order $p=0$ and $p=1$.

Such LPEs involve a kernel $K_\kappa$, taken in said study to be of von Mises-Fisher type:
\begin{equation}
   K_\kappa(x,y)=c_{\kappa,d}\exp(\kappa \langle x,y\rangle),\quad x,y\in\mathbb{S}^1,
\end{equation}
where $\kappa>0$ is a concentration parameter that needs to be selected, and $c_{\kappa,d}$ is a normalising constant, which is does not appear in the final implementation of the LPE.

Explicitly, the $p=0$, local constant LPE is given, in polar coordinates, by
\begin{equation}
    \widehat{m}_{0}(\theta)=\dfrac{\sum_{i=1}^n Y_i \exp[\kappa (X_{i,1}\cos \theta+X_{i,2}\sin \theta)]}{\sum_{i=1}^n \exp[\kappa (X_{i,1}\cos \theta+X_{i,2}\sin \theta)]},\quad \theta\in(-\pi,\pi],\quad\text{where}\;X_i=(X_{i,1},X_{i,2}).
\end{equation}

The $p=1$, local linear LPE is given by
\begin{equation}
    \widehat{m}_1(\theta)=\dfrac{\sum_{i=1}^n Y_i \exp[\kappa (X_{i,1}\cos \theta+X_{i,2}\sin \theta)]}{\sum_{i=1}^n (1+\arccos(X_{i,1}\cos \theta+X_{i,2}\sin \theta)^2)\exp[\kappa (X_{i,1}\cos \theta+X_{i,2}\sin \theta)]},\quad \theta\in(-\pi,\pi].
\end{equation}

For $p\in\{0,1\}$, the authors propose the minimiser of the leave-one-out cross validation function $\text{CV}$, for the selection of $\kappa>0$, given by
\begin{equation}
\label{eq: CV for LPE S1}
    \text{CV}(\kappa):=\dfrac{1}{n}\sum_{i=1}^n[Y_i-\widehat{m}_{p,-i}(\theta_i)]^2,
\end{equation} 
where $\widehat{m}_{p.-i}$ is the LPE using the sample with $(X_i,Y_i)$ excluded. We use the \texttt{NMinimize} function in \textit{Mathematica} in order to find the global minimum of $\text{CV}(\kappa)$, subject to $\kappa>0$. 
\vspace{0.2cm}

For each of our 30 simulations for the twelve scenarios for each regression function $R^j$, $j\in\{1,2,3,4\}$, we construct a LPE of order $p=0$ and $p=1$, selecting the value of $\kappa>0$ which minimises $\text{CV}(\kappa)$. We then compute the sample mean of the square of the $\mathbb{L}^2$-distance and compare this to the estimated MISE of our proposed KRE. These estimated MISEs are given in Table \ref{tab: LPE vs KREs S1}, alongside the sample mean of the minimisers of the $\text{CV}$ function for each estimator.

\begin{table}
\centering
    \rotatebox{90}{%
    \begin{tabular}{|c||c|c|c|}
        \multicolumn{4}{c}{\textbf{Regression function:} $R^2(\theta)$} \\
        \hline
        Estimator & $p=0$ LPE & $p=1$ LPE &KRE\\
        \hline \hline
        $n=100$, $\sigma=0.5$ & \textbf{0.737} (104.6)&  0.788 (145.8) &  0.982 (0.74) \\
        $n=100$, $\sigma=1$ &  \textbf{1.444} (52.9)&1.540 (66.8)&1.584 (0.91)\\
        $n=100$, $\sigma=2$ & 3.506 (27.3)&3.793 (34.8) & \textbf{3.175} (1.14)\\
        $n=100$, $\sigma=3$ & 6.379 (19.6)&6.718 (24.2) & \textbf{5.494} (1.32)\\
        \hline \hline
        $n=500$, $\sigma=0.5$ & \textbf{0.139} (192.2) & 0.147 (232.6) & 0.144 (0.81)\\
        $n=500$, $\sigma=1$ & 0.344 (101.5) & 0.374 (127.6)& \textbf{0.321} (1.18)\\
        $n=500$, $\sigma=2$ & 0.934 (54.2) &1.026 (69.7) & \textbf{0.737} (1.51)\\
        $n=500$, $\sigma=3$ & 1.719 (36.8)&1.906 (50.2) & \textbf{1.284} (1.69)\\
        \hline \hline
        $n=1000$, $\sigma=0.5$ & \textbf{0.069} (233.3)& 0.075 (298.5) &0.069 (1.02) \\
        $n=1000$, $\sigma=1$ &0.175 (120.4) & 0.196 (157.1)& \textbf{0.153} (1.26)\\
        $n=1000$, $\sigma=2$ &0.492 (68.1) & 0.554 (88.8)& \textbf{0.362} (1.54)\\
        $n=1000$, $\sigma=3$ & 0.929 (50.8) & 1.039 (65.2) & \textbf{0.647} (1.82)\\
        \hline
    \end{tabular}
    }%
\quad  \rotatebox{90}{%
    \begin{tabular}{|c||c|c|c|}
        \multicolumn{4}{c}{\textbf{Regression function:} $R^4(\theta)$} \\
        \hline
        Estimator & $p=0$ LPE & $p=1$ LPE &KRE\\
        \hline \hline
        $n=100$, $\sigma=0.5$ &\textbf{0.825} (148.9) & 0.843 (169.2) &  1.174 (0.58)\\
        $n=100$, $\sigma=1$ & \textbf{1.784} (70.1) & 1.827 (78.2)& 1.997 (0.66)\\
        $n=100$, $\sigma=2$ & \textbf{4.276} (29.4)&4.38 (33.7)& 4.416 (0.86)\\
        $n=100$, $\sigma=3$ & \textbf{7.552} (19.8)&7.687 (23.1)& 8.075 (1.16)\\
        \hline \hline
        $n=500$, $\sigma=0.5$ & 0.154 (275.4) & 0.158  (296.4)  &  \textbf{0.148} (0.63)\\
        $n=500$, $\sigma=1$ &  0.407 (147.5)& 0.419 (161.2)&  \textbf{0.387} (0.77)\\
        $n=500$, $\sigma=2$ & 1.127 (77.6)& 1.161 (86.3)&  \textbf{1.072} (0.96) \\
        $n=500$, $\sigma=3$ & 2.082 (51.8) & 2.12 (58.3)8& \textbf{2.020} (1.11)\\
        \hline \hline
        $n=1000$, $\sigma=0.5$ &0.081 (355.2)&0.084 (386.5)&\textbf{0.072} (0.72)\\
        $n=1000$, $\sigma=1$ &0.229 (183.3)&0.229 (201.0) & \textbf{0.196} (0.85)\\
        $n=1000$, $\sigma=2$ &0.627 (98.1) & 0.653 (109.1)&\textbf{0.580} (1.03)\\
        $n=1000$, $\sigma=3$ &1.174 (67.5) & 1.216 (75.9) & \textbf{1.093} (1.17)\\
        \hline
    \end{tabular}
    }%

    \rotatebox{90}{%
    \begin{tabular}{|c||c|c|c|}
        \multicolumn{4}{c}{\textbf{Regression function:} $R^1(\theta)$} \\
        \hline
        Estimator & $p=0$ LPE & $p=1$ LPE &KRE\\
        \hline \hline
        $n=100$, $\sigma=0.5$ & 0.374 (29.8)& 0.437 (66.0) & \textbf{0.353} (1.34)\\
        $n=100$, $\sigma=1$ & 0.795 (16.0) & 1.042 (32.8)& \textbf{0.788} (1.62)\\
        $n=100$, $\sigma=2$ & 2.059 (9.8) & 2.763 (19.5)& \textbf{1.955} (2.09) \\
        $n=100$, $\sigma=3$ &  3.961 (7.5)& 5.019 (13.4)& \textbf{3.662} (2.22)\\
        \hline \hline
        $n=500$, $\sigma=0.5$ & 0.073 (68.2) & 0.091 (125.4) & \textbf{0.068} (1.45)\\
        $n=500$, $\sigma=1$ & 0.212 (29.9)& 0.254 (70.9)& \textbf{0.165} (1.84) \\
        $n=500$, $\sigma=2$ & 0.580 (18.7)& 0.736 (39.3)& \textbf{0.465} (2.2) \\
        $n=500$, $\sigma=3$ & 1.114 (14.7) & 1.400 (27.6)& \textbf{0.926} (2.24) \\
        \hline \hline
        $n=1000$, $\sigma=0.5$ & 0.039 (84.6)& 0.051 (164.0) & \textbf{0.036} (1.51)\\
        $n=1000$, $\sigma=1$ & 0.115 (42.3)& 0.141 (90.9)& \textbf{0.079} (1.92)    \\
        $n=1000$, $\sigma=2$ &  0.287 (24.9)& 0.411 (53.2)&\textbf{0.209} (2.22)\\
        $n=1000$, $\sigma=3$ & 0.574 (18.3) &0.777 (39.8) &\textbf{0.413} (2.29) \\
        \hline
    \end{tabular}
    }%
\quad
    \rotatebox{90}{%
    \begin{tabular}{|c||c|c|c|}
        \multicolumn{4}{c}{\textbf{Regression function:} $R^3(\theta)$} \\
        \hline
        Estimator & $p=0$ LPE & $p=1$ LPE &KRE\\ 
        \hline \hline
        $n=100$, $\sigma=0.5$ & \textbf{0.614} (93.1) & 0.632 (104.5) & 0.774 (0.70)\\
        $n=100$, $\sigma=1$ & \textbf{1.388} (48.1)&1.420 (54.7) & 1.399 (0.83) \\
        $n=100$, $\sigma=2$ &  3.698 (26.1)&3.703 (29.5)&\textbf{3.450} (0.99)\\
        $n=100$, $\sigma=3$ & 6.861 (17.5)& 6.685 (20.2) & \textbf{6.525} (1.09)\\
        \hline \hline
        $n=500$, $\sigma=0.5$ &0.123 (186.3)& 0.129 (212.2) & \textbf{0.116} (0.82)\\
        $n=500$, $\sigma=1$ &  0.327 (105.1)    & 0.346 (120.7) & \textbf{0.264} (1.08) \\
        $n=500$, $\sigma=2$ & 0.918 (55.0) & 0.972 (64.2)& \textbf{0.704} (1.31)\\
        $n=500$, $\sigma=3$ &1.725 (39.4) & 1.814 (45.9)& \textbf{1.331} (1.41) \\
        \hline \hline
        $n=1000$, $\sigma=0.5$ & 0.067 (229.5) & 0.071 (266.7) & \textbf{0.059} (0.95) \\
        $n=1000$, $\sigma=1$ &0.180 (126.9) & 0.193 (148.1)& \textbf{0.143} (1.18)\\
        $n=1000$, $\sigma=2$ & 0.515 (74.0) &0.550 (85.9)&  \textbf{0.365} (1.38)\\
        $n=1000$, $\sigma=3$ & 0.972 (54.4)  & 1.030 (62.9) &  \textbf{0.715} (1.48)\\
        \hline
    \end{tabular}
    }%
    
\caption{The estimated $\text{MISE}$ of the LPE $p=0$ estimator, the LPE $p=1$ estimator, and our KRE, as estimators of the regression functions $R^j$, for $j\in\{1,2,3,4\}$, $n\in\{100,500,1000\}$, and $\sigma\in\{0.5,1,2,3\}$. The minimum estimated MISE in each row is highlighted in \textbf{bold}. In parentheses, we give the sample mean of the $\kappa>0$, in the case of the LPEs, and of the $s>0$, in the case of our KRE, which minimised their respective cross-validation function in each case. }
\label{tab: LPE vs KREs S1}
\end{table}

From the table we can see that in the case of $R^1$, our KRE is preferable across sample sizes and noise levels. In the case of $R^2$, we see that the local constant LPE has better performance for $\sigma=0.5$ and in the case $(n,\sigma)=(100,1)$, with our KRE having a lower estimated MISE in all other cases. Note that in the case $(n,\sigma)=(1000,0.5)$, the difference between the MISE of the $p=0$ LPE and our KRE is on the order of $10^{-4}$. Here, our KRE demonstrates a robustness for high levels of noise and an improvement in performance as the sample size increases. 

In the case of $R^3$, for $(n,\sigma)\in\{(100,0.5),(100,1)\}$, the $p=0$ LPE is preferable, with our KRE exhibiting better performance in the remaining cases. In the case of $R^4$, the local constant LPE has a lower error for all $n=100$ cases but our KRE is preferable for $n\in\{500,1000\}$. Both of these cases suggest that our KRE is preferable for larger sample sizes, regardless of noise levels. 
\vspace{0.2cm}

Let us further examine the difference in the behaviour of the LPEs and KRE. In Figure \ref{fig: LPE vs KRE on h, sig=1, n=1000} we plot the true regression function $R^3$ in polar coordinates, alongside the $p=0$ and $p=1$ LPEs, as well as our KRE, for a given sample of $n=1000$ points and noise level $\sigma=1$. We only plot around a local maximum (near $\theta=0$) and a local minimum (near $\theta=\pi/4$) in order to see the more subtle differences between the estimators. We can see that both LPEs provide an underestimation of the local maximum and an overestimation of the local minimum. Furthermore, neither LPE manages to capture the curvature of the extrema as accurately as the KRE.

\begin{figure}
    \centering
    \includegraphics[width=1\linewidth]{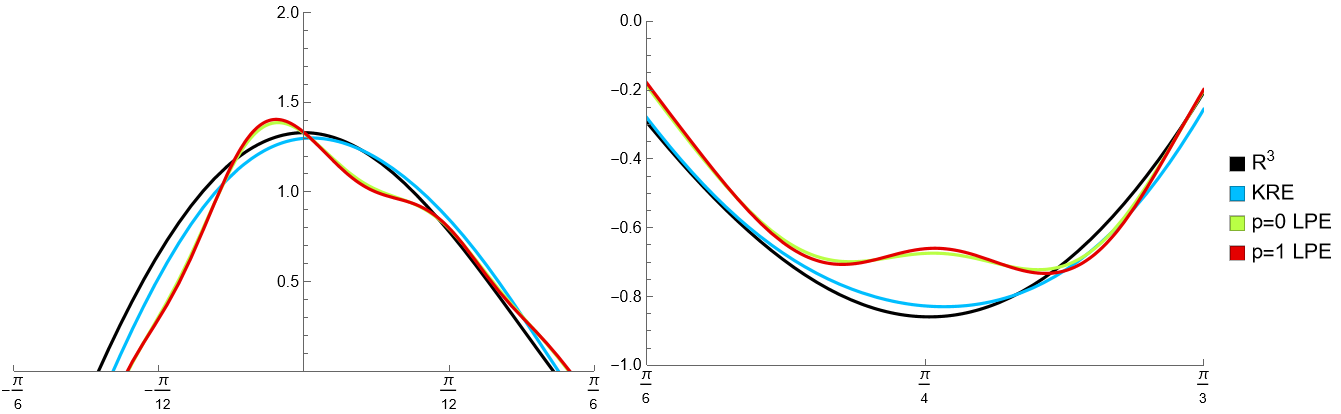}
    \caption{The regression function $R^3$ plotted around $\theta=0$ and $\theta=\pi/4$, alongside the $p=0$ LPE, the $p=1$ LPE, and the KRE, for $n=1000$ and $\sigma=1$.}
    \label{fig: LPE vs KRE on h, sig=1, n=1000}
\end{figure}

In conclusion, these results seem to suggest that the performance of the local constant polynomial estimator may be preferable in certain cases of low noise and smaller sample size, but as the noise and sample size increase, our proposed KRE appears to be a more accurate estimator. Note that all datasets discussed in Section \ref{sec: data} have sample size $n>1000$ and as such we justify the use of our KRE rather than an LPE in the following applications.

\section{Applications to Environmental Data}
\label{sec: results}
\subsection{Sea Ice Area in the Northern Hemisphere}
\label{sec: sea ice}

\begin{figure}
    \centering
    \includegraphics[width=1\linewidth]{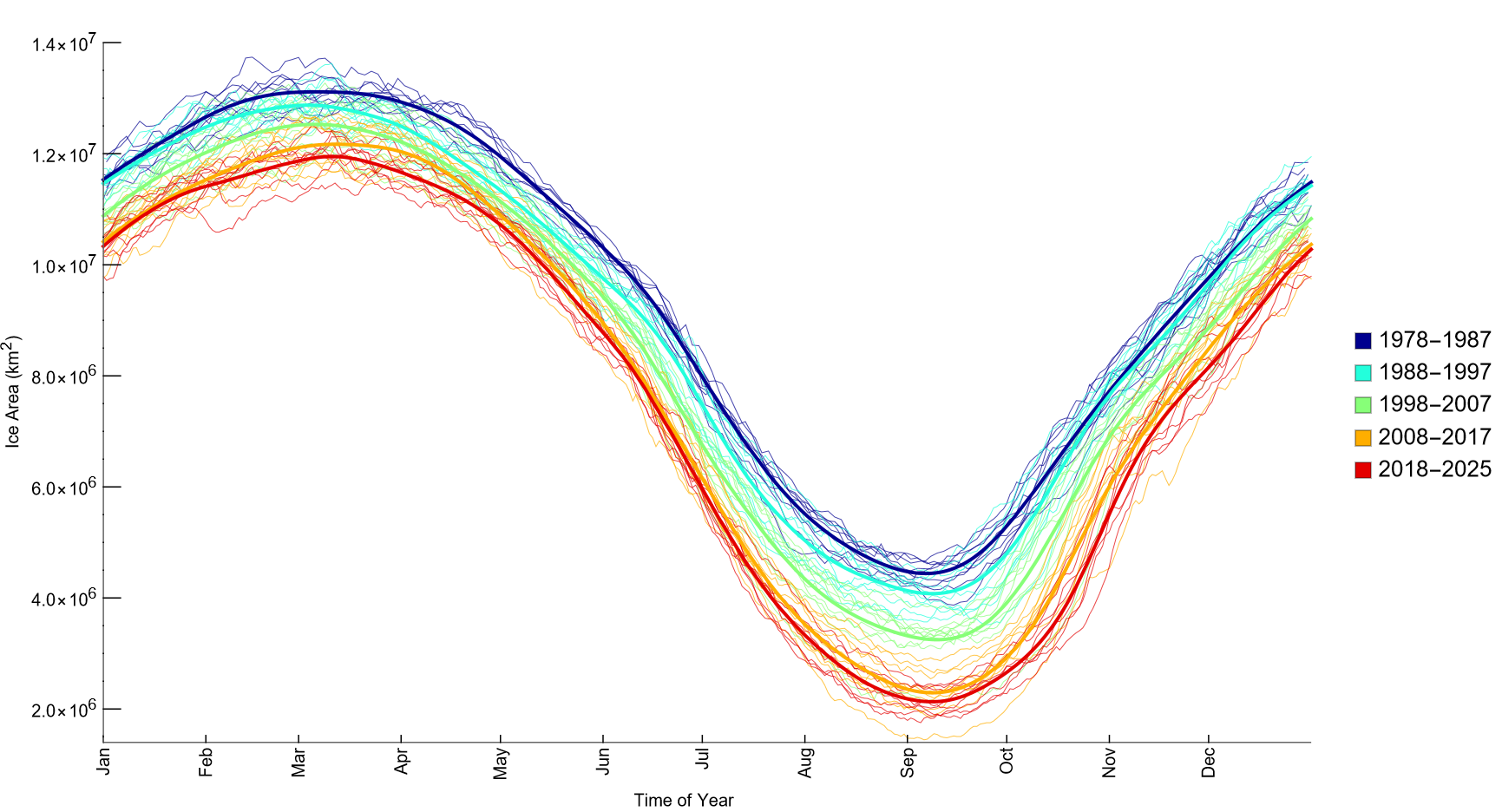}
    \caption{A plot of the observed total ice area in the Northern Hemisphere at each time of the year. The ticks on the $x$-axis indicate the first day of each month. Each of the five periods is presented in a separate colour. The thin lines on the plot represent the observed data over each period and the thick lines show the KRE based on such data.}
    \label{fig:sea ice area}
\end{figure}

For each of the five subsets (1978-1987, 1988-1997, 1998-2007, 2008-2017, 2018-2025) of the sea ice area data, we employ Algorithm \ref{alg: KRE} in order to construct a KRE. The minimiser of $\text{CV}(s)$ for each period is given in Table \ref{tab:SIA}. Such KREs are plotted in Figure \ref{fig:sea ice area} alongside the data for each year.

\begin{table}[]
    \centering
    \begin{tabular}{|c|c|c|c|c|c|c|}
    \hline
       Period  &  $\text{argmin}\,\text{CV}(s)$& $\max \widehat{R}_s$ &$\min\widehat{R}_s$&$\max\widehat{R}_s'$&$\min\widehat{R}_s'$&$\text{mean}\{e_i\}$\\
       \hline\hline
        1978-1987 & 1.25&$1.31\times10^7$&$4.44\times 10^6$&$4.82\times10^6$&$-5.58\times10^6$&$46.7$\\
        1988-1997& 1.3&$1.29\times10^7$&$4.08 \times 10^6$&$6.00\times10^6$&$-5.96\times 10^6$&$-13.0$\\
        1998-2007& 1.4&$1.25\times 10^7$&$3.25\times 10^6$&$6.23 \times10^6$&$-5.84\times 10^6$&$-3.2$\\
        2008-2017&1.35&$1.22\times 10^7$&$2.29\times 10^6$&$6.83\times10^6$&$-5.90\times10^6$&$1.3\times 10^{-10}$\\
        2018-2025&1.15&$1.20\times 10^7$&$2.13\times 10^6$&$7.88\times10^6$&$-6.28\times10^6$&$3.0\times 10^{-10}$\\
        \hline
    \end{tabular}
    \caption{The minimiser of $\text{CV}(s)$ for the KREs for each period of the sea ice area in the Northern Hemisphere dataset, the global extrema of the KREs, the extrema of the first derivative of the KREs, and the sample mean of the residuals. }
    \label{tab:SIA}
\end{table}

\begin{figure}
    \centering
    \includegraphics[width=1\linewidth]{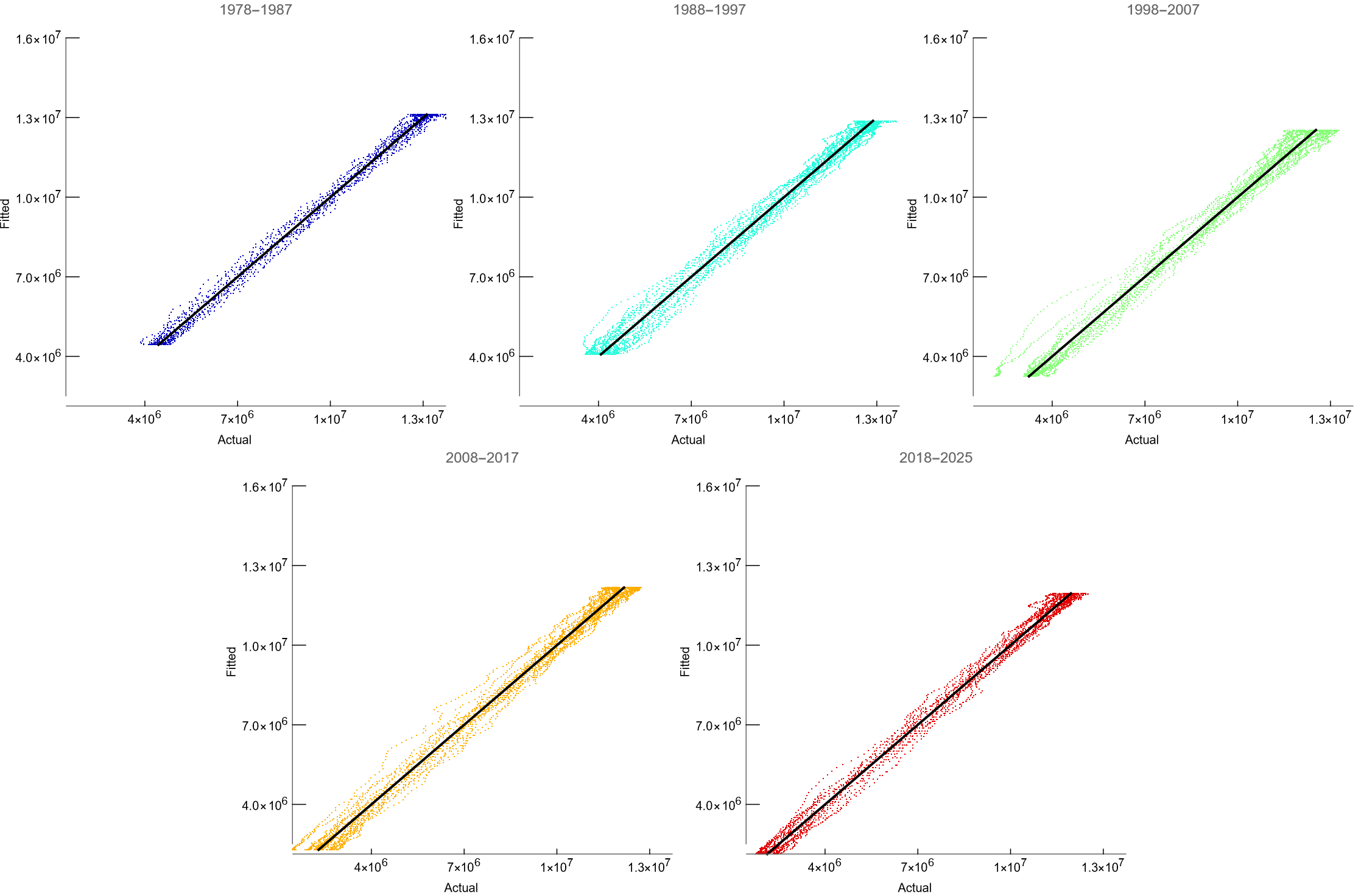}
    \caption{The actual observations $Y_i$ plotted against the fitted $\widehat{R}_s(\theta_i)$ for the sea ice area regression estimation. The black line indicates the identity function. }
    \label{fig:obs fit SIA}
\end{figure}

We assess the performance of each KRE by examining the sample residuals $e_i:=Y_i-\widehat{R}_s(X_i)$, $i=1,\dots,n$. Figure \ref{fig:obs fit SIA} presents the actual versus fitted values, which appear well scattered around the identity line. The sample mean residuals for each period are given in Table \ref{tab:SIA}. Relative to the scale of the response variable ($\simeq 10^6$), these means are close to zero. Overall, the results suggest a good fit of our regression model.

Within Figure \ref{fig:sea ice area}, we can see that from decade to decade the sea ice area is decreasing. Indeed, each consecutive KRE lies entirely below the KRE for the previous period. We can see that the decrease in sea ice area is more dramatic near the global minimum, in September. This is in line with current scientific findings on the levels of sea ice in the Arctic \cite{DriversSeaIce}. This drastic decrease in summer ice levels across decades is indicative of the future ice-free summer predicted by \cite{DK21}.

\begin{figure}
    \centering
    \includegraphics[width=1\linewidth]{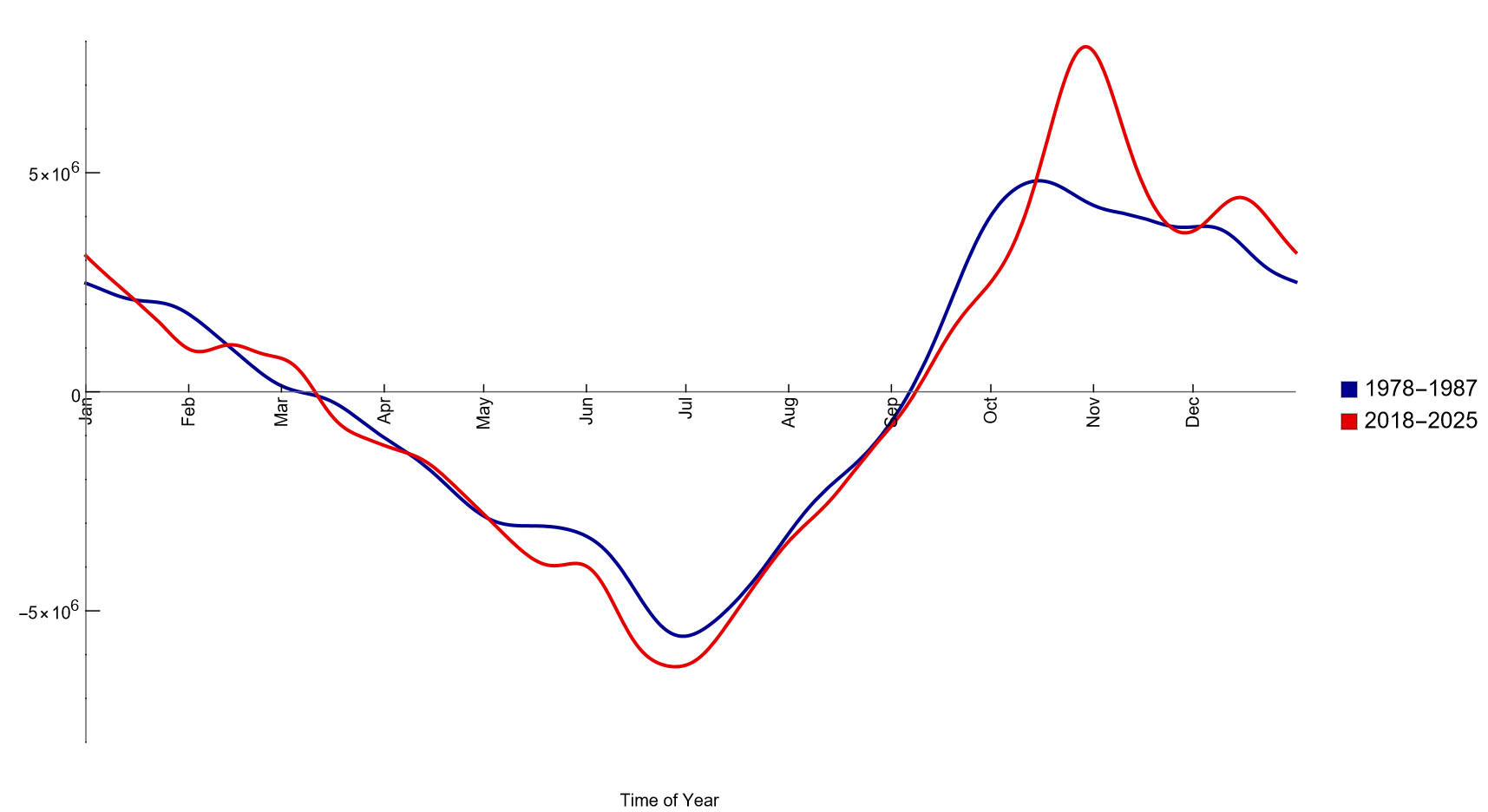}
    \caption{The first derivative of the the KRE for the period 1978-1987 and for 2018-2025. The ticks on the $x$-axis indicate the first day of each month.}
    \label{fig:deriv of SIA}
\end{figure}

An important advantage of our methodology, compared with the use of average curves as in \cite{SL16}, is that it allows us to study not only the estimated regression functions and their differences across historical periods, but also their derivatives. This provides additional information about the rates and timing of seasonal changes. Each derivative has two zeros in its domain, corresponding to a single season of freezing and a single season of thawing, indicating one period of increasing and one period of decreasing sea ice area during the annual cycle. In Table \ref{tab:SIA} we give the extrema of the first derivative of each KRE. The maximum of the first derivative increases across historical periods, while its minimum becomes more negative generally. This suggests that, in addition to the decline in sea ice area, the maximum estimated rates of both refreezing and retreat have increased in magnitude. These changes in seasonal dynamics are illustrated in Figure \ref{fig:deriv of SIA}, where we compare the first derivatives for 1978--1987 and 2018--2025.

\subsection{Daily Mean Temperature in Stockholm}
\label{sec: stockholm}

\begin{figure}
    \centering
    \includegraphics[width=1\linewidth]{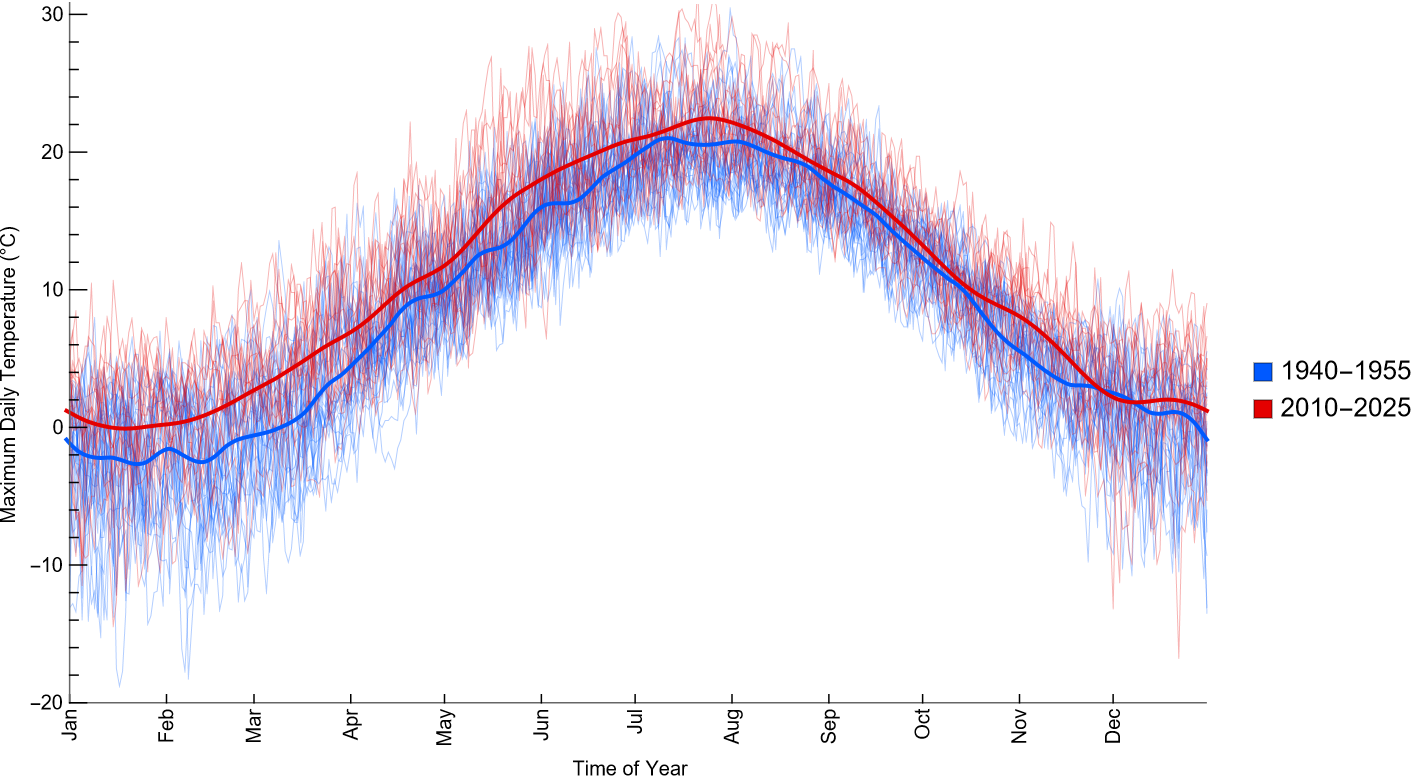}
    \caption{The thin lines represent observed daily mean temperature and the thick lines give the KREs. The period 1940-1955 is shown in blue and 2010-2025 is in red. The ticks on the time axis indicate the first day of each month.}
    \label{fig:temp stockholm}
\end{figure}

For both periods (1940-1955 and 2010-2025) of daily mean temperature data in Stockholm, we construct a KRE as in Algorithm \ref{alg: KRE}. The values of $s>0$ minimising $\text{CV}(s)$ are given in Table \ref{tab:stockholm}. We denote the 1940-1955 KRE by $\widehat{R}_{\text{past}}$, and the 2010-2025 KRE by $\widehat{R}_\text{present}$. Both KREs are plotted in Figure \ref{fig:temp stockholm} alongside the observed curves for each year.

\begin{table}[]
    \centering
    \begin{tabular}{|c|c|c|c|c|c|c|}
    \hline
       Period & $\text{argmin}\,\text{CV}(s)$ & $\max \widehat{R}_s$ & $\min\widehat{R}_s$ & $\max D_{\delta}\widehat{R}_s(\theta)$ & $\min D_{\delta}\widehat{R}_s(\theta)$ & $\text{mean}\{e_i\}$\\
       \hline\hline
        1940-1955 & 0.95 & $20.99$ & $-2.68$ & $10.61$ & $-10.46$ & $-3.4\times10^{-16}$\\
        2010-2025 & 0.9 & $22.46$ & $-0.09$ & $9.72$ & $-10.47$ & $3.3\times10^{-16}$\\
    \hline
    \end{tabular}
    \caption{The minimiser of $\text{CV}(s)$ for the daily mean temperature in Stockholm in each period, the estimated global extrema of each KRE, and the sample mean of the residuals.}
    \label{tab:stockholm}
\end{table}

As in the sea ice application, an examination of the residuals suggests a satisfactory fit of the KREs to the observed data. The sample mean residuals for both periods are given in Table \ref{tab:stockholm} and are close to zero.

\begin{figure}
    \centering
    \includegraphics[width=1\linewidth]{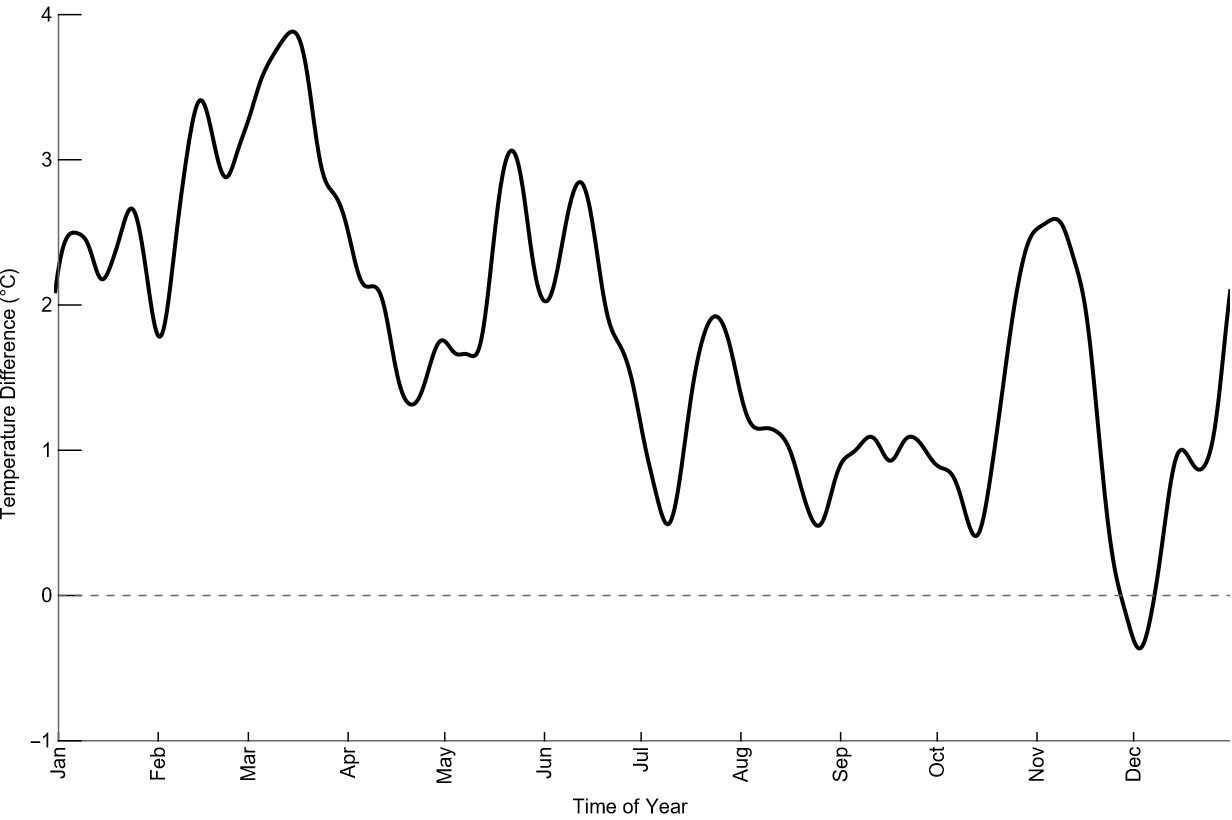}
    \caption{A plot of $\widehat{R}_\text{present}-\widehat{R}_\text{past}$ over the course of the year for the Stockholm temperature data.}
    \label{fig:stock diff}
\end{figure}

The KREs both capture the seasonal weather patterns of Sweden with warm summers and cold winters. Note that the two KREs are noticeably different with the present KRE indicating higher average daily temperatures across most of the year. Indeed, in Figure \ref{fig:stock diff} we plot $\widehat{R}_\text{present}-\widehat{R}_\text{past}$ over the course of a year. 

For the vast majority ($\approx96.8\%$) of the year the 2010-2025 KRE has a higher daily mean temperature; as much as $\approx 3.88^\circ\text{C}$ warmer at time. Indeed, the only period where the KRE $\widehat{R}_\text{present}$ gives \textit{lower} daily mean temperatures than $\widehat{R}_\text{past}$ is approximately the period from the $26^\text{th}$ of November until the $7^\text{th}$ of December. In this period, the lowest difference in temperature is $\approx 0.36^\circ\text{C}$.

The derivatives of the KREs give us the approximate minimum and maximum of each KRE. These values are given in Table \ref{tab:stockholm}. We can see a noticeable increase in both the global minima and maxima between the periods. Moreover numerical integration tells us that the average value of the difference of the two KREs over the course of a year is
\begin{equation}
    \dfrac{1}{2\pi}\int_{-\pi}^\pi \big(\widehat{R}_\text{present}(\theta)-\widehat{R}_\text{past}(\theta)\big) \,d\theta\approx1.78^\circ \text{C},
\end{equation}
which indicates a significant increase in the daily mean temperature across the year. This is above than the target of the 2015  Paris Agreement to pursue efforts to keep long-term increases in global average temperatures to less than $1.5^\circ$C compared to `pre-industrial levels'. We firstly remark that this is only a local result and not necessarily indicative of global trends. Moreover, the recent work of \cite{Thorne26} lays out the nuances needed when interpreting such results as they relate to conclusions about the extent of human-induced warming. Nonetheless, this is a noteworthy result indicating the increase in average temperatures across the course of the year in Stockholm, since 1940.

\begin{figure}
     \centering
     \includegraphics[width=1\linewidth]{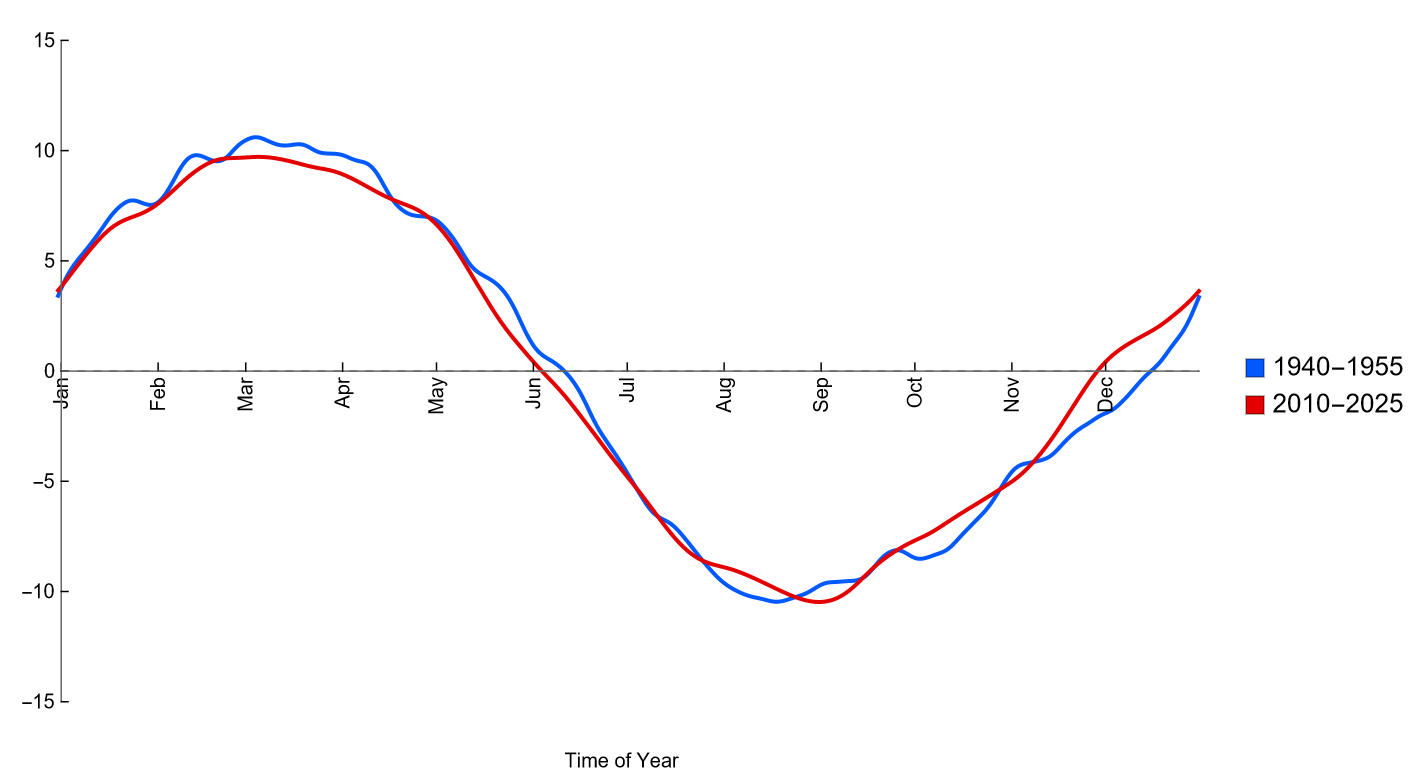}
     \caption{The discrete derivative of the Stockholm temperature KREs with $\delta=2\pi/4$ over the course of the year.}
     \label{fig:disc deriv stock}
 \end{figure}

To further investigate the \textit{seasonal dynamics} of temperature, we examine changes in the fitted regression functions over fixed time intervals. Unlike the instantaneous rates described by ordinary derivatives, finite differences allow us to study temperature changes over a chosen seasonal time scale. We define
the following `discrete derivatives':
 \begin{equation}
     D_{\delta}f(x):=\dfrac{f(x+\delta)-f(x)}{\delta},
 \end{equation}
where $\delta>0$ is a fixed increment. If $f$ is periodic, then so is $D_\delta f$. We take $\delta=2\pi/4$, corresponding to approximately one quarter of a year, to examine seasonal temperature changes. Figure \ref{fig:disc deriv stock} presents the resulting finite differences for the two historical periods. This analysis compliments the comparison of fitted temperature levels by providing information about changes in their seasonal evolution. Interestingly, the finite differences exhibit similar seasonal patterns and extrema for the two historical periods, despite the increase in fitted temperature levels.
 
From these results, we have evidence to support the claim that past data of the daily mean temperature of Stockholm has a different distribution to more present data. 

\section{Conclusions}
\label{sec: conclusion}

In this paper, we developed a general kernel-based approach to nonparametric regression, using kernel density estimators as building blocks. Focusing on the unit circle, we constructed explicit and computationally tractable kernel regression estimators for periodic data, together with a leave-one-out cross-validation procedure for the data-driven selection of the smoothness parameter. The estimators are straightforward to implement, and the publicly available code (see Section \ref{sec: data avail}) facilitates their use in other applications. An extensive simulation study demonstrated the effectiveness of the proposed estimators, including their competitive performance compared with the literature.

We applied the proposed methodology to two important environmental indicators. The analysis of Arctic sea ice area revealed a systematic decline across five historical periods between 1978 and 2025, while the Stockholm temperature application indicated higher fitted temperatures during 2010--2025 than during 1940--1955 across approximately $97\%$ of the year, with an estimated annual average increase of $1.78^\circ$C. 

Beyond these changes in the levels of the environmental indicators, we investigated their \textit{seasonal dynamics} through derivatives and finite differences of the fitted regression functions. The results suggest increasing magnitudes of the maximum seasonal rates of Arctic sea ice retreat and refreezing, whereas the seasonal temperature changes in Stockholm exhibit broadly similar patterns across the two historical periods, despite the increase in fitted temperature levels.

Our methodology provides a basis for further developments in nonparametric regression on more general metric spaces, whenever suitable kernel density estimators are available. Further directions include investigating the theoretical properties and convergence rates of the proposed regression estimators, as well as developing statistical inference procedures for comparing fitted regression functions and their seasonal dynamics across historical periods.

\section{Data \& Code Availability}
\label{sec: data avail}
All data used in Section \ref{sec: data} is freely available; the Sea Ice Concentrations data from the National Snow and Ice Data Center is found at \url{https://nsidc.org/data/nsidc-0051/versions/2}. The daily mean temperature in Stockholm data is accessible via the Open-Meteo Historical Weather API, \url{https://open-meteo.com/en/docs/historical-weather-api}. 

\texttt{R} code for the implementation of Algorithm \ref{alg: KRE} is available at \url{https://github.com/APPercival/KREs-on-S1}.

\end{document}